\documentclass{amsart}
\title{Unicity of types and local Jacquet--Langlands correspondence}
\author{Yuki Yamamoto}
\usepackage{mathrsfs}
\usepackage{amssymb}
\usepackage[all]{xy}
\newtheorem{defn}{Definition}[section]

\newtheorem{thm}[defn]{Theorem}
\newtheorem{prop}[defn]{Proposition}

\newtheorem{eg}[defn]{Example}

\DeclareMathOperator{\GL}{GL}

\DeclareMathOperator{\M}{M}

\DeclareMathOperator{\cInd}{c--Ind}

\DeclareMathOperator{\Gal}{Gal}

\DeclareMathOperator{\JL}{JL}
\DeclareMathOperator{\LJ}{LJ}
\DeclareMathOperator{\Irr}{Irr}
\DeclareMathOperator{\MT}{MT}
\newcommand{\Afra}{\mathfrak{A}}
\newcommand{\Bfra}{\mathfrak{B}}
\newcommand{\Bcal}{\mathcal{B}}

\newcommand{\Mcal}{\mathcal{M}}

\newcommand{\ofra}{\mathfrak{o}}

\newcommand{\Z}{\mathbb{Z}}
\newcommand{\Q}{\mathbb{Q}}

\newcommand{\C}{\mathbb{C}}

\begin{document}
\maketitle

\begin{abstract}
Let $F$ be a non-archimedean local field.  
For any irreducible representation $\pi$ of an inner form $G'=\GL_{m}(D)$ of $G=\GL_{N}(F)$, there exists an irredubile representation of a maximal compact open subgroup in $G'$ which is also a type for $\pi$.  
Then we can consider the problem whether these types are unique or not in some sense.  
If such types for $\pi$ are unique, we say $\pi$ has the strong unicity property of types.  
On the other hand, there exists a correspondence connecting irreducible representations of $G'$ and $G$, called the Jacquet--Langland correspondence.  
In this paper, we study the ralation between the strong unicity of types and the Jacquet--Langlands correspondence.  
\end{abstract}

\tableofcontents

\section{Introduction}

Let $F$ be a non-archimedean local field.  
Let $m \in \Z_{>0}$, and let $D$ be a central division $F$-algebra.  
In the following, we consider smooth representations of locally profinite groups over the field of complex numbers.  

To study representation theory of $G'=\GL_{m}(D)$, the Bernstein decomposition is important.  
We recall it.  
We denote by $\Mcal(G')$ the category of $G'$-representations.  
Let $\Irr(G')$ be the set of isomorphism classes of irreducible representations of $G'$ over $\C$.  
Then there exists an equivalent relation $\sim$ on $\Irr(G')$, called ``the same inertial support''.  
We put $\Bcal(G'):=\Irr(G')/\sim$, and we denote by $\mathfrak{s}$ the quotient map $\Irr(G') \to \Bcal(G')$.  
Then, for $\mathfrak{s}_{0} \in \Bcal(G')$, we can consider the full subcategory $\Mcal(G')_{\mathfrak{s}_{0}}$ of $\Mcal(G')$, whose objects are $G'$-representations $\pi$ such that any irreducible subquotient of $\pi$ is an element in $\mathfrak{s}_{0}$.  
In the setting, it is known that we have the decomposition of the category
\[
	\Mcal(G') = \prod_{\mathfrak{s}_{0} \in \Bcal(G')} \Mcal(G')_{\mathfrak{s}_{0}}, 
\]
called the Bernstein decomposition.  
Moreover, this decomposition is indecomposable.  
In this sense, the Bernstein decomposition is a concept to classify irreducible representations of $G'$.  

As a concept connecting with the Bernstein decomposition, there exists a class of representations of some compact open subgroups in $G'$, called types.  
Let $\mathfrak{s}_{0} \in \Bcal(G')$.  
We say a pair $(J, \lambda)$ consisting a compact open subgroup $J$ in $G$ and an irreducible $J$-representation $\lambda$ is an $\mathfrak{s}_{0}$-type if 
\[
	\pi \in \mathfrak{s}_{0} \text{ if and only if } \lambda \subset \pi|_{J}
\]
for any $\pi \in \Irr(G')$.  
Then, types are useful to examine whether an irreducible representation $\pi$ is an element in $\mathfrak{s}_{0}$ or not by studying the restriction of $\pi$ to some compact open subgroup in $G'$.  

For $G'=\GL_{m}(D)$, for any $\pi \in \mathcal{A}_{m}(D)$ there exists an $\mathfrak{s}(\pi)$-type.  
For example, S\'echerre \cite{S1}, \cite{S2}, \cite{S3}, \cite{SS} constructed types for $G'$, called simple types.  
It is a generalization of the construction by Bushnell--Kutzko \cite{BK} for $G=\GL_{N}(F)$ case.  

Here, we can also consider another question whether for an irreducible representation $\pi'$ of $G'$, types for $\pi$ are unique or not.  
This question is first considered in \cite[Appendix A]{BM} to formulate the Breuil--M\'ezard conjecture for $G=\GL_{2}(\Q_{p})$ case.  
For any irreducible representation $\pi$ of $\GL_{2}(F)$, Henniart calculated the number of irreducible representations $\rho$ of $K=\GL_{2}(\ofra_{F})$ such that $(K, \rho)$ is also an $\mathfrak{s}(\pi)$-type, where $\ofra_{F}$ is the ring of integers in $F$.  

To generalize this problem to more general reductive group case, we introduce archetypes, defined by Latham \cite[Definition 4.3]{La}.  
We denote by $\MT(\mathfrak{s}_{0})$ the set of $\mathfrak{s}_{0}$-types $(K, \rho)$ such that $K$ is a maximal compact open subgroup in $G'$.  
The group $G$ acts $\MT(\mathfrak{s}_{0})$ by conjugation.  
We call these orbits $\mathfrak{s}_{0}$-archetypes.  

When $G'=\GL_{N}(F)$ and $\pi$ is irreducible and supercuspidal, Pa\v{s}k\={u}nas showed in \cite{Pas} that there exists a unique $\mathfrak{s}(\pi)$-archetype.  
On the other hand, it is shown in \cite[A.1.5(3)]{BM} that if $F$ satisfies some assumption, there exists a non-supercuspidal representation $\pi \in \mathcal{A}_{2}(F)$ such that $\mathfrak{s}(\pi)$-archetypes are not unique.  
Let $G'=\GL_{m}(D)$.  
If $\pi \in \Irr(G')$ is supercuspidal and `unramified' in some sense, the author \cite{YY1} showed that there exists a unique $\mathfrak{s}(\pi)$-archetype.  
On the other hand, the author also constructed a supercuspidal representation $\pi$ of $\GL_{m}(D)$ for some $m$ and $D$ such that $\mathfrak{s}(\pi)$-archetypes are not unique.  

\begin{defn}[{{\cite[\S 1]{LN}}}]
Let $\pi$ be an irreducible representation of $G$.  
We say $\pi$ has \textit{the strong unicity property of types} (SUP) if there exists a unique $\mathfrak{s}(\pi)$-archetype.  
\end{defn}

Here, we denote by $\mathcal{A}_{m}(D)$ the set of isomorphism classes of irreducible, essentially square-integrable representations of $\GL_{m}(D)$.  
We put $N=m \cdot (\dim_{F}D)^{1/2}$ and $G=\GL_{N}(F)$.  
Then there exists a correspondence 
\[
	\JL_{m,D} \colon \mathcal{A}_{m}(D) \overset{\sim}{\to} \mathcal{A}_{N}(F)
\]
connecting representations of $G'$ and $G$, called the Jacquet--Langlands correspondence.  
It is known that $\JL_{m,D}$ satisfies some `nice' properties.  
We also put $\LJ_{m,D}:=(\JL_{m,D})^{-1}$.  

Then we have a natural question:  ``Does the Jacquet--Langlands correspondence preserve the strong unicity of types?''
The goal of this paper is to study the behavior of the strong unicity property of types by the Jacquet--Langlands correspondence.  

Our main theorem is the following.  

\begin{thm}[Theorem \ref{MainThm}]
\begin{enumerate}
\item For some $F$, $m$ and $D$, the map $\JL_{m,D}$ does not preserve the strong unicity property of types.  
\item For some $F$, $m'$ and $D'$, the map $\LJ_{m',D'}$ does not preserve the strong unicity property of types.  
\end{enumerate}
\end{thm}

\bigbreak
\noindent{\bfseries Acknowledgments}\quad
I am deeply grateful to my supervisor Naoki Imai for giving helpful comments to improve the draft.  
I also thank Yuanqing Cai for his insightful suggestions.  
I was supported by the FMSP program at Graduate School of Mathematical Sciences, the University of Tokyo.  
I am also supported by JSPS KAKENHI Grant Number JP21J13751.  
\bigbreak

\section{Preliminaries}

\subsection{Notation}

For a division $F$-algebra $D$, we denote by $\ofra_{D}$ the ring of integers in $D$.  

Let $G$ be a locally profinite group, and let $H$ be a closed subgroup in $G$.  
For $g \in G$, we put ${}^{g}H=gHg^{-1}$.  
For a representation $\rho$ of $H$, we define a ${}^{g}H$-representation ${}^{g} \rho$ by ${}^{g}\rho(ghg^{-1})=\rho(h)$ for $h \in H$.  

\subsection{Simple type}

We recall simple types for inner forms of $\GL_{N}$, which are types for essentially square-integrable representations, from \cite{S1}, \cite{S2}, \cite{S3}, \cite{SS}.  

Simple types are constructed from information of a 4-tuple $[\Afra, n, 0, \beta]$, called a simple stratum.  
Simple strata $[\Afra, n, 0, \beta]$ satisfies some conditions, for example, $\beta$ is an element in $\M_{m}(D)$ such that $F[\beta]$ is an extension field of $F$.  
From $[\Afra, n, 0, \beta]$, we can define compact open subgroups $J=J(\beta, \Afra)$, $J^{1}=J^{1}(\beta, \Afra)$ and $H^{1}=H^{1}(\beta, \Afra)$ in $G'=\GL_{m}(D)$.  
To construct simple types, we need two constructions of irreducible $J$-representations:  
\begin{enumerate}
\item We can define a finite set $\mathscr{C}(\beta, \Afra)$ of `simple characters' from $[\Afra, n, 0, \beta]$.  
When we take $\theta \in \mathscr{C}(\beta, \Afra)$, we can define a unique irreducible $J^{1}$-representation $\eta_{\theta}$ up to isomorphism.  
Moreover, we can take a `nice' extension $\kappa$ of $\eta_{\theta}$ to $J$, called a $\beta$-extension.  
\item Suppose $J/J^{1}$ is isomorphic to $\GL_{l}(\mathbf{f})^{\times j}$, where $\mathbf{f}$ is a finite field.  
When $\eta$ is an irreducible cuspidal representation of $\GL_{l}(\mathbf{f})$, we can take the inflation $\sigma$ of $\eta^{\otimes j}$ to $J$.  
\end{enumerate}
When we take $\kappa$ and $\sigma$ as above, we put $\lambda = \kappa \otimes \sigma$, and then $(J, \lambda)$ is a simple type.  
If $(J, \lambda)$ is a type for some supercuspidal representation, we say $(J, \lambda)$ is a maximal simple type.  

From two distinct simple strata $[\Afra, n, 0, \beta]$ and $[\Afra', n', 0, \beta']$, we may construct the same simple type.  
In this case, it is shown in \cite[Theorem 9.4]{BSS} that $[F[\beta]:F]=[F[\beta']:F]$.  

\subsection{Parametric degree}
For $\pi \in \mathcal{A}_{m}(D)$, we can define a positive integer $\delta(\pi)$ as in \cite[\S 2]{BH}, called the parametric degree of $\pi$.  

We recall the definition of $\delta(\pi)$ for some supercuspidal representation $\pi$.  
Let $(J, \lambda)$ be a maximal simple type which is also an $\mathfrak{s}(\pi)$-type.  
Suppose $(J, \lambda)$ is associated with a simple stratum $[\Afra, n, 0, \beta]$.  
Then there exist two irreducible representations $\kappa$ and $\sigma$ of $J$ such that 
\begin{itemize}
\item $\kappa$ is constructed from a `simple character' associated with $[\Afra, n, 0, \beta]$, 
\item $\sigma$ is the inflation of some irreducible representation of a product of general linear groups over some finite field, and
\item $\lambda = \kappa \otimes \sigma$.  
\end{itemize}
In this case, we can define a positive integer $\delta_{0}(\sigma_{\Bfra})$ as in \cite[2.4, 2.6]{BH} and we put $\delta_{0}(\lambda):=\delta_{0}(\sigma_{\Bfra})[F[\beta]:F]$.  

\begin{prop}[{{\cite[Proposition 2.7]{BH}}}]
\begin{enumerate}
\item The integer $\delta_{0}(\lambda)$ is independent of the choice of simple strata and representations $\sigma$ and $\kappa$ such that $\lambda \cong \kappa \otimes \sigma$.  
\item If $(J, \lambda)$ and $(J', \lambda')$ are simple types which are also $\mathfrak{s}(\pi)$-types for some irreducible representation $\pi$, then $\delta_{0}(\lambda)=\delta_{0}(\lambda')$.  
\end{enumerate}
\end{prop}

Then, for $\pi' \in \mathcal{A}_{m}(D)$ which is supercuspidal, we can define a positive integer $\delta(\pi')$ by $\delta(\pi')=\delta_{0}(\lambda)$ for a maximal simple type $(J, \lambda)$ which is also an $\mathfrak{s}(\pi)$-type.  

We can generalize the parametric degree $\delta(\pi)$ to any $\pi \in \mathcal{A}_{m}(D)$.  
The integer $\delta(\pi)$ satisfies the following properties.  

\begin{prop}[{{\cite[2.7.3, 2.8 Corollary 1]{BH}}}]
\label{Prop_for_delta}
\begin{enumerate}
\item For any $\pi \in \mathcal{A}_{m}(D)$, the integer $\delta(\pi)$ divides $N=m \cdot (\dim_{F}D)^{1/2}$.  
\label{dividingN}
\item For $\pi \in \mathcal{A}_{N}(F)$, the representation $\pi$ is supercuspidal if and only if $\delta(\pi)=N$.  
\label{delta_and_sc}
\item For any $m$ and $D$, the map $\JL_{m,D}$ preserves the parametric degree, that is, we have $\delta(\JL_{m,D}(\pi'))=\delta(\pi')$ for any $\pi' \in \mathcal{A}_{m}(D)$.  
\label{delta_and_JL}
\end{enumerate}
\end{prop}

\section{Main examples}

\begin{eg}
We consider $N=2$ and $G'=D^{\times}$, where $D$ is a quaternion algebra of $F$.  
Suppose the order of the residual field of $F$ is $2$.  
Then, by \cite[A 1.5(3)]{BM} there exists an essentially square-integrable representation $\pi$ of $G=\GL_{2}(F)$ and irreducible representations $\rho_{1}$ and $\rho_{2}$ of $K=\GL_{2}(\ofra_{F})$ such that $\rho_{1} \ncong \rho_{2}$ and $(K, \rho_{i})$ is an $\mathfrak{s}(\pi)$-type for $i=1,2$.  

To show $\pi$ does not have SUP, suppose $\rho_{1}$ and $\rho_{2}$ are $G$-conjugate.  
Then there exists $g \in G$ such that $K=gKg^{-1}$ and $\rho_{2} \cong {}^{g}\rho_{1}$.  
Since the normalizer of $K$ is $F^{\times}K$, there exist $t \in F^{\times}$ and $k \in K$ such that $g=tk$.  
Here, we have ${}^{k} \rho_{1} \cong \rho_{1}$ since $\rho_{1}$ is a $K$-representation and $k \in K$.  
Moreover, we also have ${}^{t} \rho_{1} \cong \rho_{1}$ since $t$ is an element in the center of $G$.  
Therefore $\rho_{2} \cong {}^{g} \rho_{1} = {}^{tk} \rho_{1} \cong {}^{t} \rho_{1} \cong \rho_{1}$, which is a contradiction.  

On the other hand, we put $\pi':=\LJ_{1,D}(\pi)$.  
Since $\pi'$ is an irreducible representation of $D^{\times}$, which is the multiplicative group of a central division $F$-algebra, the representation $\pi'$ has SUP.  
That is, although $\pi'$ has SUP, the image $\pi=\JL_{1,D}(\pi')$ does not have SUP.  
Therefore, $\JL_{1,D}$ does not preserve the strong unicity of types.  
\end{eg}

\begin{eg}
We consider $N=4$ and $G'=\GL_{2}(D)$, where $D$ is a quaternion algebra of $F$.  
Let $E_{1}/F$ be a ramified quadratic extension of fields, and let $E/E_{1}$ be an unramified quadratic extension of fields.  
We consider an embedding $E \hookrightarrow \M_{2}(D)$ as in \cite[Example 7.1]{YY1}.  
We can take an element $\beta \in E$ such that $E=F[\beta]$ and a 4-tuple $[\Afra, n, 0, \beta]$ is a simple stratum, where $\Afra = \M_{2}(\ofra_{D})$ and $n \in \Z_{>0}$ is defined from $\Afra$ and $\beta$.  
Let $(J, \lambda)$ be a maximal simple type associated with the simple stratum $[\Afra, n, 0, \beta]$, and let $(\tilde{J}, \Lambda)$ be a maximal extension of $(J, \lambda)$ in $G$.  
We put $\pi':=\cInd_{\tilde{J}}^{G} \Lambda$.  
Then $\pi'$ is irreducible and supercuspidal, and $\pi'$ does not have SUP by \cite[\S 7]{YY1}.  

We put $\pi:=\JL_{2,D}(\pi')$.  
We will show $\pi$ has SUP.  
Since any irreducible supercuspidal representation of $\GL_{N}(F)$ has SUP by \cite{Pas}, it is enough to show $\pi$ is supercuspidal.  
Then we consider the parametric degree $\delta(\pi)$ of $\pi$.  
To show $\pi$ is supercuspidal, it suffices to show $\delta(\pi)=N=4$ by Proposition $\ref{Prop_for_delta} (\ref{delta_and_sc})$. 
Since $\delta(\pi)=\delta(\JL(\pi'))=\delta(\pi')$ by Proposition $\ref{Prop_for_delta} (\ref{delta_and_JL})$, it is enough to show $\delta(\pi')=4$.  
Here, we have $[F[\beta]:F]=[E:F]=4$.  
Since $\delta(\pi')=\delta_{0}(\lambda)$ is a multiple of $[F[\beta]:F]$ by the definition of $\delta_{0}(\lambda)$, we have $\delta(\pi') \in 4\Z$.  
On the other hand, the integer $\delta(\pi')$ divides $4$ by Proposition $\ref{Prop_for_delta} (\ref{dividingN})$.  
Then $\delta(\pi')=4$ and $\pi$ is supercuspidal.  
That is, although $\pi$ has SUP, the image $\pi'=\LJ_{2,D}(\pi)$ does not have SUP.  
Therefore, $\LJ_{2,D}$ does not preserve the strong unicity of types.  
\end{eg}

As a consequence of these examples, we obtain the following theorem.  

\begin{thm}
\label{MainThm}
\begin{enumerate}
\item Suppose the order of the residual field of $F$ is $2$.  
Let $D$ be a quaternion algebra of $F$.  
Then $\JL_{1,D}$ does not preserve the strong unicity property of types, that is, there exists $\pi' \in \mathcal{A}_{1}(D)$ such that $\pi'$ has SUP and $\JL_{1,D}(\pi') \in \mathcal{A}_{2}(F)$ does not have SUP.  
\item Let $D$ be a quaternion algebra of $F$.  
Then $\LJ_{2,D}$ does not preserve the strong unicity property of types, that is, there exists $\pi \in \mathcal{A}_{4}(F)$ such that $\pi$ has SUP and $\LJ_{2,D}(\pi) \in \mathcal{A}_{2}(D)$ does not have SUP.  
\end{enumerate}
\end{thm}

\end{document}